\documentclass{ifacconf}
\usepackage{graphicx}      
\usepackage{natbib}        
\usepackage{graphics} 
\usepackage{epsfig} 
\usepackage{amsmath} 
\usepackage{amssymb}  
\usepackage{amsfonts}  
\usepackage{subfig}
\usepackage{epstopdf}
\usepackage{enumerate}
\usepackage{graphicx} 

\usepackage{color}
\newcommand{\tcb}{\textcolor{blue}}

\newcommand*{\QEDB}{\hfill\ensuremath{\square}}%

\newtheorem{theorem}{Theorem}

\newtheorem{definition}{Definition}
\newtheorem{assumption}{Assumption}
\newtheorem{remark}{Remark}

\begin{document}
\begin{frontmatter}

\title{Singularly Perturbed Stochastic Hybrid Systems:\\ Stability and Recurrence via Composite Nonsmooth Foster Functions\thanksref{footnoteinfo}}

\thanks[footnoteinfo]{Research supported in part by AFOSR grant FA9550-22-1-0211 and NSF grant ECCS CAREER 2305756.\\E-mail address: poveda@ucsd.edu}

\author{Jorge I. Poveda}

\address{Department of Electrical and Computer Engineering,\\
University of California, San Diego, CA 92093, USA}

\begin{abstract}               
We introduce new sufficient conditions for verifying stability and recurrence properties in singularly perturbed \emph{stochastic hybrid dynamical systems}. Specifically, we focus on hybrid systems with deterministic continuous-time dynamics that exhibit multiple time scales and are modeled by constrained differential inclusions, as well as discrete-time dynamics modeled by constrained difference inclusions with random inputs. By assuming regularity and causality of the dynamics and their solutions, respectively, we propose a suitable class of composite nonsmooth Lagrange-Foster and Lyapunov-Foster functions that can certify stability and recurrence using simpler functions related to the slow and fast dynamics of the system. We establish the stability properties with respect to compact sets, while the recurrence properties are studied only for open sets. 
\end{abstract}

\begin{keyword}
Stochastic hybrid systems, Stochastic Stability, Lyapunov Functions.
\end{keyword}

\end{frontmatter}

\section{Introduction}

When modeling dynamical systems as ordinary differential equations (ODEs), singular perturbation tools are often used to analyze systems in which some of the states, known as ``fast states", rapidly converge to a quasi-steady state manifold parameterized by the ``slow states". In such cases, the stability properties of the original system can typically be predicted by examining the stability properties of a reduced system that ignores the fast dynamics. For ODEs, these tools were introduced around the 1940's in (\tcb{\cite{krylov1947introduction}}), and have been further developed in the references (\tcb{\cite{Kokotovic_SP_Book,Saberi,TeelMoreauNesic2003}}), etc. 

On the other hand, singular perturbation tools for deterministic hybrid systems that combine continuous-time and discrete-time dynamics were studied in \tcb{\cite{Sanfelice:11,Wang:12_Automatica,SP2020Teel,OnlinePovedaHybrid}}, and \tcb{\cite{RajebSP}}. These tools have enabled the development of new hybrid controllers and algorithms for systems that exhibit multiple time scales in the flows, including hybrid and non-smooth extremum seeking control (\tcb{\cite{PoTe17Auto,PovedaKrsticFXT}}), distributed optimization algorithms under attacks (\tcb{\cite{wang2020stability}}), and accelerated distributed learning algorithms for network games (\tcb{\cite{TAC21Momentum_Nash}}). However, even though some of these deterministic hybrid control systems have overcome the limitations of smooth feedback-based algorithms, there are still many engineering, biological, financial, and mathematical problems that involve \emph{stochasticity}. Examples include algorithms designed to achieve robust global convergence to the global minimizer of a non-convex cost (\tcb{\cite{baradaran2018stochastic}}), solving model-free optimization problems with stochastic perturbations (\tcb{\cite{Ramirez_SourceSeeking}}), coordinating multi-agent systems on manifolds (\tcb{\cite{HartmanBook}}), achieving steady-state Nash equilibrium seeking in games with plants in the loop via stochastic perturbations (\tcb{\cite{PoTe16CDC}}), and synchronizing oscillators with binary phase update rules (\tcb{\cite{JavedAutomatica20}}). In all these applications, the closed-loop system is modeled as a \emph{stochastic hybrid dynamical system} (SHDS), as described in \tcb{\cite{teel2015stochastic}} and in \tcb{\cite{Teel:14_AutomaticaSurvey}}.

While singular perturbation tools have been widely used for ODEs and deterministic hybrid systems, their application to singularly perturbed stochastic hybrid dynamical systems (SP-SHDS) has remained mostly unexplored. Existing research has mainly focused on certain classes of Markov linear systems \tcb{\cite{HaddadsingularlyPerturbed,Hinfinitesingularly}}, finite-time horizon problems \tcb{\cite{FilarTAC}}, and stochastic ODEs. This gap in the literature motivates our paper, which introduces new tools for the stability analysis of a class of SP-SHDS that combine deterministic differential inclusions and stochastic difference inclusions. Our main contribution is the study of two properties of SHDS that have not been considered before in the context of singular perturbations: \emph{Uniform Global Asymptotic Stability in Probability} (UGASp), and \emph{Uniform Global Recurrence} (UGR). The former defines a stochastic extension of traditional uniform global asymptotic stability (UGAS) properties studied for deterministic hybrid and non-hybrid systems (\tcb{\cite{bookHDS}}), while the latter is studied mostly for stochastic systems \tcb{\cite{MeybBook}}. To certify these properties, we study nonsmooth Lyapunov-based tools based on composite regular Lagrange and Foster functions constructed from simpler functions related to the slow and fast dynamics of the system. Our approach leverages the stability tools introduced in the seminal work of \tcb{\cite{Teel_ANRC}} for general well-posed SHDS, and it is the first Lyapunov-based stability characterization for SP-SHDS with set-valued dynamics, as well as the first study of the property of recurrence in such systems. The certificates provided in the paper build upon suitable composite constructions that expand on previous results in the literature established for ODEs in \tcb{\cite{Saberi,NaiduSP}} and for deterministic hybrid systems in \tcb{\cite{Wang2020SP}}. The tools outlined in the paper can facilitate the creation and assessment of novel forms of multi-time scale control, optimization, and estimation algorithms that involve both hybrid dynamics and stochasticity in the jumps. 
\section{PRELIMINARIES}
\label{sec:preli}

\vspace{-0.2cm}
\subsection{Notation and Basics on Nonsmooth Analysis}

\vspace{-0.2cm}
We denote the set of (non negative) real numbers by $(\mathbb{R}_{\geq 0})$ $\mathbb{R}$.
The set of (nonnegative) integers is denoted by $(\mathbb{Z}_{\geq 0})$ $\mathbb{Z}$. Given a 
set $\mathcal{A} \subset \mathbb{R}^n$ and a vector $z \in \mathbb{R}^n$, we define $|z|_\mathcal{A} := \inf_{y \in \mathcal{A}}|z-y|$, and we use $|\cdot|$ to denote the standard Euclidean norm. We use $\overline{\mathcal{A}}$ to denote the closure of $\mathcal{A}$, and $
r\mathbb{B}^\circ$ to denote the open ball (in the Euclidean norm) of appropriate dimension centered around the origin and with radius $r>0$. For ease of notation, given two vectors $u,v \in \mathbb{R}^{n}$, we write $(u,v)$ for $(u^{T},v^{T})^{T}$. A function $f:\mathbb{R}^n\to\mathbb{R}$ is said to be: a)  $C^k$ if its $k^{th}$ derivative is continuous; and b) radially unbounded if $f(x)\to \infty$ whenever $|x|\to\infty$. A function $\alpha:\mathbb{R}_{\geq0}\to\mathbb{R}_{\geq0}$ is said to be: a) of class $\mathcal{G}_{\infty}$ if it is continuous, non-decreasing, and unbounded; b) of class $\mathcal{K}_{\infty}$ if it is zero at zero, continuous, strictly increasing, and unbounded. It is said to be of class $\mathcal{P}_s\mathcal{D}(\mathcal{A})$ if it is positive semidefinite with respect to $\mathcal{A}$, and of class $\mathcal{P}\mathcal{D}(\mathcal{A})$ when it is positive definite with respect to $\mathcal{A}$. When $\mathcal{A}=\{0\}$, we simply use $\mathcal{P}_s\mathcal{D}$ and $\mathcal{P}\mathcal{D}$. We use $\overline{\text{co}}(\mathcal{A})$ to denote the closure of the convex hull of the set $\mathcal{A}$, and $\mathbb{I}_{\mathcal{A}}:\mathbb{R}^n\to\{0,1\}$ to denote the standard indicator function. We use $\mathbf{B}(\mathbb{R}^m)$ to denote the Borel $\sigma$-field, and $K\subset\mathbb{R}^m$ is said to be measurable if $K\in\mathbf{B}(\mathbb{R}^m)$.

Let $f:\mathbb{R}^n\to\mathbb{R}$ be a locally Lipschitz function, and let $\mathcal{Z}$ be the set of points where $\nabla f$ is not defined, which is of measure zero due to Rademacher's Theorem. The Clarke generalized gradient of $f$ at $y\in\text{dom}~f$, is the set
\begin{equation*}
\partial f(y):=\text{co}\{v\in\mathbb{R}^n:\exists y_k\to y,~y_k\notin\mathcal{Z},~\lim_{k\to\infty}\nabla f(x_k)=v\}.
\end{equation*}
We use $\partial_{x_1} f(x_1,x_2)$ and $\partial_{x_2} f(x_1,x_2)$ to denote the partial Clarke gradients. The function $f$ is said to be regular at $y$ if, for every $u\in\mathbb{R}^n$, the directional derivative $f'(y;u):=\lim_{s\to 0^+}\frac{f(y+su)-f(y)}{s}$ exists, and $f'(x;u)=\max\{v^\top u:v\in \partial f(x)\}$, for all $u\in\mathbb{R}^n$. Typical examples of locally Lipschitz regular functions include $\mathcal{C}^1$ and convex functions. A set-valued mapping $F:\mathbb{R}^m\rightrightarrows\mathbb{R}^n$ is outer semi-continuous (OSC) if for each $(x_i,y_i)\to(x,y)\in\mathbb{R}^m\times\mathbb{R}^n$ satisfying $y_i\in M(x_i)$ for all $i\in\mathbb{Z}_{\geq0}$, we have $y\in M(x)$. A mapping $F$ is locally bounded (LB) if, for each bounded set $K$, $M(K):=\bigcup_{x\in K}M(x)$ is bounded.  Given a set $\mathcal{X} \subset \mathbb{R}^m$, the mapping $M$ is OSC and LB relative to $\mathcal{X}$ if the set-valued mapping from $\mathbb{R}^m$ to $\mathbb{R}^n$ defined by $M$ for $x \in \mathcal{X}$,  and by $\varnothing$ for $x \notin \mathcal{X}$, is OSC and LB at each $x \in \mathcal{X}$. The graph of $F$ is the set $\text{graph}(F) := \{(x, y) \in  \mathbb{R}^m \times \mathbb{R}^n: y \in F(x)\}$. Given a measurable space $(\Omega, \mathcal{F})$, a set-valued map $F: \Omega \rightrightarrows \mathbb{R}^n$ is said to be \textit{$\mathcal{F}$-measurable}, if for each open set $\mathcal{O} \subset \mathbb{R}^n$, the set $F^{-1}(\mathcal{O}) :=
\{\omega \in \Omega : F(\omega) \cap \mathcal{O} 	= \varnothing\}\in \mathcal{F}$.  

\subsection{Basic Notions on Stochastic Hybrid Dynamical Systems}
\label{prelimSHDS}
In this paper, we consider SHDS of the form
\begin{subequations}\label{SHDS1}
	\begin{align}
	&y\in C,~~~~~~~~~\dot{y}\in F(y),\label{SHDS_flows0}\\
	&y\in D,~~~~~~y^+\in G(y,v^+),~~~v\sim \mu(\cdot),~~\label{SHDS_jumps0}
	\end{align}
\end{subequations}
where $F:\mathbb{R}^n\rightrightarrows\mathbb{R}^n$ is called the flow map, $G:\mathbb{R}^n\times\mathbb{R}^m\rightrightarrows\mathbb{R}^n$ is called the jump map, $C$ is the flow set, $D$ is the jump set, and $v^+$ is a place holder for a sequence $\{\bf{v}_k\}_{k=1}^{\infty}$ of independent, identically distributed (i.i.d.) random variables $\bf{v}_k:\Omega\to\mathbb{R}^m$, $k\in\mathbb{N}$, defined on a probability space $(\Omega,\mathcal{F},\mathbb{P})$. Thus, ${\bf v_k}^{-1}(F):=\{\omega\in\Omega: \mathbf{v}_k(\omega)\in F\}\in\mathcal{F}$ for all $F\in\mathbf{B}(\mathbb{R})^m$, and $\mu:\mathbf{B}(\mathbb{R}^m)\to[0,1]$ is defined as $\mu(F):=\mathbb{P}\{\omega\in\Omega:{\bf v}_k(\omega)\in F\}$. 

For the purpose of completeness, we review the concept of solution to \eqref{SHDS1}, presented in \tcb{\cite{Teel_ANRC}}. Random solutions to SHDS \eqref{SHDS1} are functions of $\omega\in\Omega$, denoted ${\bf y}(\omega)$. To formally define these mappings, for $\ell\in\mathbb{Z}_{\geq1}$, let $\mathcal{F}_\ell$ denote the collection of sets $\{\omega\in\Omega:({\bf v}_1(\omega),{\bf v}_2(\omega),\ldots,{\bf {v}}_\ell(\omega))\in F\}$, $F\in\mathbf{B}(\mathbb{R}^m)^\ell)$, which are the sub-$\sigma$-fields of $\mathcal{F}$ that form the minimal filtration of ${\bf v}=\{{\bf v}_\ell \}_{\ell=1}^{\infty}$, which is the smallest $\sigma$-algebra on $(\Omega,\mathcal{F})$ that contains the pre-images of $\mathbf{B}(\mathbb{R}^m)$-measurable subsets on $\mathbb{R}^m$ for times up to $\ell$. A stochastic hybrid arc is a mapping ${\bf y}$ from $\Omega$ to the set of hybrid arcs \tcb{\cite[Ch. 2]{bookHDS}}, such that the set-valued mapping from $\Omega$ to $\mathbb{R}^{n+2}$, given by  $\omega\mapsto \text{graph}({\bf y}(\omega)):=\big\{(t,j,z):\tilde{y}={\bf y}(\omega), (t,j)\in\text{dom}(\tilde{y}),z=\tilde{y}(t,j)\big\}$, is $\mathcal{F}$-measurable with closed-values. Let $\text{graph}({\bf y}(\omega))_{\leq \ell}:=\text{graph}({\bf y} (\omega))\cap (\mathbb{R}_{\geq0}\times\{0,1,\ldots,\ell\}\times\mathbb{R}^n)$. An $\{\mathcal{F}_\ell\}_{\ell=0}^{\infty}$ adapted stochastic hybrid arc is a stochastic hybrid arc ${\bf y}$ such that the mapping $\omega\mapsto \text{graph}({\bf y}(\omega))_{\leq \ell}$ is $\mathcal{F}_\ell$ measurable for each $\ell \in\mathbb{N}$. An adapted stochastic hybrid arc ${\bf y}(\omega)$, or simply $\mathbf{y}_\omega$,  is a solution to the SHDS \eqref{SHDS1} starting from $y_0$, denoted ${\bf y}_\omega\in \mathcal{S}_r(y_0)$, if: (1) $\mathbf{y}_\omega(0,0)=y_0$; (2) if $(t_1,j),(t_2,j)\in\text{dom}(\mathbf{y})$ with $t_1<t_2$, then for almost all $t\in[t_1,t_2]$, $\mathbf{y}_{\omega}(t,j)\in C$ and $\dot{\mathbf{y}}_\omega(t,j)\in F(\mathbf{y}_\omega(t,j))$; (3) if $(t,j),(t,j+1)\in\text{dom}(\mathbf{y}_\omega)$, then $\mathbf{y}_\omega(t,j)\in D$ and $\mathbf{y}_\omega(t,j+1)\in G(\mathbf{y}_\omega(t,j),\mathbf{v}_{j+1}(\omega))$. A random solution ${\bf y}_\omega$ is said to be: a) almost surely {\em non-trivial} if its hybrid time domain contains at least two points almost surely; and b) almost surely  {\em  complete} if for almost every sample path $\omega\in \Omega$ the hybrid arc ${\bf y_\omega}$ has an unbounded time domain. 
%
\section{A Class of Singularly Perturbed SHDS}
\label{section3}
We focus our attention on a sub-class of SHDS \eqref{SHDS1}, given by

\vspace{-0.8cm}
\begin{subequations}\label{SPSHDS1}
\begin{align}
&(x,z)\in C,~~\dot{x}\in F_x(x,z),~~~\varepsilon\dot{z}\in F_z(x,z)
\label{flow_dynamics1}\\
&(x,z)\in D,~~(x^+,z^+)\in G(x,z,v^+),~~v\sim \mu,\label{jump_dynamics1}
\end{align}
\end{subequations}
where $\varepsilon\in\mathbb{R}_{>0}$ is a small parameter, $x\in\mathbb{R}^{n_1}$ is the ``slow'' state, $z\in\mathbb{R}^{n_2}$ is the ``fast'' state, $F_x:\mathbb{R}^{n_1}\times\mathbb{R}^{n_2}\rightrightarrows\mathbb{R}^{n_1}$, $F_z:\mathbb{R}^{n_1}\times\mathbb{R}^{n_2}\rightrightarrows\mathbb{R}^{n_2}$ are set-valued mappings, and $C,D\subset\mathbb{R}^{n_1}\times\mathbb{R}^{n_2}$ define the flow set and the jump sets, respectively. To simplify our presentation, we will consider sets of the form $C:=C_x\times C_z$ and $D:=D_x\times D_z$, where $C_x,D_x\subset\mathbb{R}^{n_1}$ and $C_z,D_z\subset\mathbb{R}^{n_2}$. We will also use $$F_{\varepsilon}(x,z):=F_x(x,z)\times \frac{1}{\varepsilon}F_z(x,z),$$ to denote the overall flow map, $y=(x,z)$ to denote the overall state, and $\mathcal{V}:=\bigcup_{\omega\in\Omega,i\in\mathbb{Z}_{\geq0}} \mathbf{v}_{i+1}(\omega)$ to denote the set of all possible values that $v$ can take.

In system \eqref{SPSHDS1}, $v$ is a place-holder for a sequence of i.i.d random variables with probability measure $\mu$. Since $\varepsilon=0$ causes a singularity in the flow dynamics \eqref{flow_dynamics1}, system \eqref{SPSHDS1} is a \emph{singularly perturbed stochastic hybrid dynamical system} (SP-SHDS). This model covers different types of systems previously studied in the literature. In particular: 
\begin{enumerate}[(a)]
\item When $G$ is independent of the random input $v$, system \eqref{SPSHDS1} recovers the model studied in \tcb{\cite{OnlinePovedaHybrid,SP2020Teel}} in the context of deterministic hybrid dynamical systems with set-valued flows, and in \tcb{\cite{RajebSP}} for linear hybrid systems. 
\item If, in addition to (a), the maps $F_x$,$F_z$ are singled-valued and continuous, and the flow set and the jump set components of $z$ are compact, system \eqref{SPSHDS1} recovers the models studied in \tcb{\cite{Sanfelice:11,PoTe17Auto}}. 
\item If, in addition to (a) and (b), $D=\emptyset$, system \eqref{SPSHDS1} recovers the models studied in \tcb{\cite{SP_Survey2012,HistorySingularPerturbations}} in the context of singularly perturbed differential inclusions.
\item If, in addition to (a), (b), and (c), $F_x$ and $F_z$ are locally Lipschitz, then system \eqref{SPSHDS1} recovers the models studied in \tcb{\cite{Saberi}} for ODEs.
\item Finally, if in addition to (a)-(d) the dynamics \eqref{SPSHDS1} are linear, then system \eqref{SPSHDS1} recovers the models studied in \tcb{\cite[Ch. 1-6]{Kokotovic_SP_Book}}.
\end{enumerate}
The incorporation of stochastic dynamics into singularly perturbed hybrid dynamical systems, particularly those having set-valued jump maps, can induce a complete spectrum of behaviors that have not been fully studied before in the literature of singular perturbations. Of particular interest to us are \emph{causal} behaviors that are consistent with the definition of random solutions given in Section \ref{prelimSHDS}. 

\vspace{0.1cm}
 \begin{definition}\label{definitionbasic1}
The SP-SHDS \eqref{SPSHDS1} is said to satisfy the \textsl{Basic Conditions} if \tcb{\cite{Teel_ANRC}}: (a) $C$ and $D$ are closed,  $C\subset \text{dom}(F_{\varepsilon})$, and $D\subset\text{dom}(G)$. (b) The set-valued mapping $F_{\varepsilon}$ is OSC, LB, and convex-valued for any $\varepsilon>0$. (c) The set-valued map $G$ is LB and the mapping $v\mapsto \text{graph}(G(\cdot,\cdot,v)):=\{(s_0,s_1,s_2)\in\mathbb{R}^{n_1}\times\mathbb{R}^{n_2}\times\mathbb{R}^{n_1+n_2}:s_2\in G(s_0,s_1,v)\}$ is measurable with closed values. \QEDB
\end{definition}
In this paper, we will only work with SHDS that satisfy the Basic Conditions: 

\noindent 
\textsl{Standing Assumption.} The SP-SHDS \eqref{SPSHDS1}, \tcb{and its reduced dynamics \eqref{SPSHDSreduced}}, satisfy the Basic Conditions. \QEDB

Our goal is to predict the stability properties of \eqref{SPSHDS1} based on the properties of a simpler \emph{reduced} SHDS characterized by the steady state condition of the fast dynamics in \eqref{flow_dynamics1}, also known as the \emph{boundary-layer} dynamics.
\subsection{Boundary Layer Dynamics}
The boundary layer dynamics of \eqref{SPSHDS1} ignore the jumps, and are given by the system
\begin{equation}\label{BLD}
(x,z)\in C,~~\dot{z}\in F_{z}(x,z),~~\dot{x}=0.
\end{equation}
Note that in \eqref{BLD}, the state $x$ is fixed. For these dynamics, we will assume the existence of a quasi-steady-state map $\mathcal{M}$, possible set-valued, which, for each fixed $x$, generates a set that is uniformly attractive and stable for system \eqref{BLD}.
\begin{assumption}\label{manifold}
There exists an OSC and LB set-valued mapping $\mathcal{M}:\mathbb{R}^{n_1}\rightrightarrows\mathbb{R}^{n_2}$, such that for all $x\in C_x$, $\mathcal{M}(x)\neq\emptyset$ and $\mathcal{M}(x)\subset C_z$. \QEDB
\end{assumption}
To characterize the stability properties of the boundary layer dynamics with respect to $\mathcal{M}$, we will use the following assumption. 

\begin{assumption}\label{Assumption2}
There exists a locally Lipschitz and regular function $W:\mathbb{R}^{n_1}\times\mathbb{R}^{n_2}\to\mathbb{R}_{\geq0}$, $\alpha_1,\alpha_2\in\mathcal{G}_{\infty}$, a continuous function $\varphi_z:\mathbb{R}_{\geq0}\to\mathbb{R}_{\geq0}$, and $k_z\in\mathbb{R}_{>0}$ such that:
\begin{enumerate}[(a)]
\item For all $y\in C\cup D\cup G(D\times\mathcal{V})$, we have:
\begin{equation}\label{Kinftyboundsfast}
\alpha_1\left(|z|_{\mathcal{M}(x)}\right)\leq W(y)\leq \alpha_2\left(|z|_{\mathcal{M}(x)}\right).
\end{equation}
\item For all $y\in C$ and all $f_z\in F_{z}(y)$, we have:
\begin{equation}\label{boundaryflowsLyapunov}
\max_{\substack{\nu\in \partial_z W(y)}} \langle \nu,f_z\rangle  \leq -k_z \varphi_z^2\left(|z|_{\mathcal{M}(x)}\right).
\end{equation}
\end{enumerate}
\end{assumption}
\begin{remark}
For each fixed $x$, the conditions of Assumption \ref{Assumption2} essentially establish uniform Lagrange stability of the set $\mathcal{M}(x)$ for the differential inclusion \eqref{BLD}. If, additionally, $\alpha_1,\alpha_2\in\mathcal{K}_{\infty}$ and $\varphi_z$ is positive definite, they imply Uniform Global Asymptotic Stability \tcb{\cite[Ch.3]{bookHDS}}.
\end{remark}
\begin{remark}
When $\mathcal{M}$ is a single-valued continuous function, Assumptions \ref{manifold}-\ref{Assumption2} recover the standard assumptions considered in the literature of ODEs, where $\mathcal{M}$ is usually a quasi-steady state manifold, see \tcb{\cite{Saberi,Kokotovic_SP_Book,TeelMoreauNesic2003}}. 
\end{remark}
In some applications, the stochastic jumps of the fast state $z$ might also converge towards $x$ in a probabilistic sense. In that particular situation, the following Assumption will be considered.
\begin{assumption}\label{assumptionjumpsfast}
There exists $c_z>0$ and $\rho_z\in\mathcal{P}_s\mathcal{D}$ such that the function $W$ of Assumption \ref{Assumption2} satisfies:
\begin{equation}
\int_{\mathbb{R}^m}\sup_{\substack{g\in G(y,v)}}~W(g)\mu(dv)\leq W(y)-c_z\rho_z(|z|_{\mathcal{M}(x)}).
\end{equation}
for all $y\in D$. \QEDB
\end{assumption}
In words, Assumption \ref{assumptionjumpsfast} asks that the worst case value of the function $W$ evaluated during jumps of the \emph{original} SHDS \eqref{SPSHDS1} does not increase in expectation.

\subsection{Reduced Stochastic Hybrid Dynamics}
Using Assumption \ref{manifold}, we define the \emph{reduced} SHDS associated with \eqref{SPSHDS1}, given by
\begin{subequations}\label{SPSHDSreduced}
\begin{align}
&x\in C_x,~~~~\dot{x}\in \tilde{F}(x),\label{flow_dynamics1r}\\
&x\in D_x,~~x^+\in \tilde{G}(x,v^+),~~v\sim \mu,\label{jump_dynamics1r}
\end{align}
\end{subequations}
where $\tilde{F}$ and $\tilde{G}$ are
\begin{subequations}
\begin{align}
&\tilde{F}(x):=\overline{\text{co}}\left\{\tilde{f}_x\in\mathbb{R}^{n_1}:\tilde{f}_x\in F_x(x,z),z\in \mathcal{M}(x)\right\},\\
&\tilde{G}(x,v):=\left\{s\in\mathbb{R}^{n_1}:(s,l)\in G(x,z,v),z\in D_z\right\}.
\end{align}
\end{subequations}
Note that this reduced hybrid system, which ignores the dynamics of $z$, is also stochastic via \eqref{jump_dynamics1r}. 

To capture the dynamic properties of \eqref{SPSHDSreduced}, we will first consider the following assumption on the flows.

\begin{assumption}\label{Assumption3}
There exists a locally Lipschitz and regular function $V:\mathbb{R}^{n_1}\to\mathbb{R}_{\geq0}$, functions $\alpha_3,\alpha_4\in \mathcal{G}_{\infty}$, a continuous function $\varphi:\mathbb{R}^{n_1}\to\mathbb{R}_{\geq0}$, and constants $k_x\in\mathbb{R}_{>0},\mu_F\in\mathbb{R}_{\geq0}$, such that:
\begin{enumerate}[(a)]
\item For all $x\in C_x\cup D_x \cup \tilde{G}(D_x\times\mathcal{V})$, we have:
\begin{equation}\label{slowKinfty}
\alpha_3(|x|_{\mathcal{A}})\leq V(x)\leq \alpha_4(|x|_{\mathcal{A}}),
\end{equation}
where $\mathcal{A}\subset\mathbb{R}^{n_1}$ is a compact set.
\item For all $x\in C_x$ and all $\tilde{f}_x\in \tilde{F}(x)$, we have:
\begin{equation}\label{reducedflowsLyapunov}
\max_{\substack{\nu\in \partial V(x)}}\langle \nu,\tilde{f}_x\rangle \leq -k_x \varphi_x^2(x)+\mu_F \mathbb{I}_{\mathcal{O}}(x),
\end{equation}
where $\mathcal{O}\subset\mathbb{R}^{n_1}$ is an open and bounded set. \QEDB
\end{enumerate}
\end{assumption}
\begin{remark}
As in Assumption \ref{Assumption2}, the two conditions of Assumption \ref{Assumption3} essentially guarantee uniform Lagrange stability of the set $\mathcal{A}$ for the differential inclusion \eqref{flow_dynamics1r} whenever $\mu_F=0$.  When $\mu_F>0$, condition \eqref{reducedflowsLyapunov} will allow us to study recurrence properties of the open set $\mathcal{O}$.
\end{remark}

For some applications, the reduced SHDS \eqref{SPSHDSreduced} might also exhibit suitable stability properties during the jumps. In this particular case, the following assumption will be used in conjunction with Assumption \ref{Assumption3}.
\begin{assumption}\label{assumptionjumpsreduced}
There exists a continuous function $\rho_x:\mathbb{R}^{n_1}\to\mathbb{R}_{\geq0}$, $c_x\in\mathbb{R}_{>0}$, $\mu_J\in\mathbb{R}_{\geq0}$, such that  
\begin{equation}\label{jumpsreducedstochastic}
\int_{\mathbb{R}^m}\sup_{\substack{\tilde{g}\in \tilde{G}(x,v)}}~V(\tilde{g})\mu(dv)\leq V(x)-c_x\rho_x(x)+\mu_J\mathbb{I}_{\mathcal{O}}(x),
\end{equation}
for all $x\in D_x$. \QEDB
\end{assumption}
If Assumptions \ref{Assumption3} and \ref{assumptionjumpsreduced} hold with $\mathcal{A}=\overline{\mathcal{O}}$, 
then we can directly conclude uniform global recurrence (see Def. 3 in Section 4.2) of the set $\mathcal{O}$ for the \emph{reduced} SHDS \eqref{SPSHDSreduced} via \tcb{\cite[Thm. 4.4]{Teel_ANRC}}. If, additionally, $\alpha_3,\alpha_4\in\mathcal{K}_{\infty}$, $\rho_x\in\mathcal{P}\mathcal{D}(\mathcal{A})$, and $\mu_J=0$, then we can conclude UGASp (see Def. 2 of Section 4.1) of the set $\mathcal{A}$ for the \emph{reduced} SDHS via \tcb{\cite[Thm. 4.5]{Teel_ANRC}}. Whether or not these properties will be preserved, in some sense, by the original SP-SHDS \eqref{SPSHDS1}, will depend on the value of $\varepsilon$ and some additional ``interconnection'' conditions.

\vspace{-0.2cm}
\subsection{Interconnection Conditions}

\vspace{-0.2cm}
We will consider two different types of interconnection conditions: one related to the flows \eqref{flow_dynamics1}, and the other one related to the jumps \eqref{jump_dynamics1}. We shall not necessarily use both conditions simultaneously. 

Below, in Assumptions \ref{assumption6}-\ref{interconnection2} the functions $V$, $W$, $\varphi_x$, and $\varphi_z$ are the same from Assumptions \ref{Assumption2} and \ref{Assumption3}.

\begin{assumption}\label{assumption6}
There exist $k_1,k_2,k_3\in\mathbb{R}_{>0}$, such that:
\begin{enumerate}[(a)]
\item For all $y\in C$, and for all $f_x\in F_x(x,z)$, we have:
\begin{equation}\label{coupled1Lyapu1}
\max_{\substack{\nu\in \partial_x W(y)}} \langle \nu, f_x \rangle \leq k_1\varphi_z(|z|_{\mathcal{M}(x)})\varphi_x(x)+k_2\varphi_z^2(|z|_{\mathcal{M}(x)}).
\end{equation}
\item For all $y\in C$, and for all $f_x\in F_x(x,z)$, there exists $\tilde{f}_x\in \tilde{F}_x(x)$ such that:
\begin{equation}\label{coupled1Lyapu2}
\max_{\substack{\nu\in \partial V(x)}}\langle \nu,f_x -\tilde{f}_x \rangle \leq k_3 \varphi_z(|z|_{\mathcal{M}(x)})\varphi_x(x).
\end{equation}
\end{enumerate}
\vspace{-0.3cm}
\QEDB 
\end{assumption}
\begin{assumption}\label{jumpcoupledcondition}
There exist  $k_4\in\mathbb{R}_{>0}$ and a continuous function $\rho_4:\mathbb{R}^{n_1}\to\mathbb{R}_{\geq0}$, such that
\begin{equation}\label{jumpcondition_expectation}
\int_{\mathbb{R}^m}\sup_{\substack{g\in G(y,v)}}~W(g)\mu(dv)\leq W(y)+k_4\rho_4(x), 
\end{equation}
for all $y\in D$. \QEDB
\end{assumption}
\begin{assumption}\label{interconnection2}
There exist  $k_5\in\mathbb{R}_{>0}$, and  $\rho_5\in\mathcal{P}_s\mathcal{D}$, such that for all $y\in D$:
\begin{equation}
\int_{\mathbb{R}^m}\sup_{\substack{\tilde{g}_x\in \tilde{G}(x,v)}}~V(g_x)\mu(dv)\leq V(x)+k_5\rho_5(|z|_{\mathcal{M}(x)}).
\end{equation} 
\vspace{-0.2cm}\QEDB
\end{assumption}
When $D=G=\emptyset$, $\mu_F=0$, $\mathcal{A}=\{0\}$, $\varphi_x,\varphi_z\in\mathcal{P}\mathcal{D}$, and all the vector fields are singled-valued and locally Lipschitz, the conditions of Assumptions \ref{Assumption2}, \ref{assumptionjumpsfast}, and \ref{assumption6} essentially recover the quadratic-type characterization presented in \tcb{\cite{Saberi}} for Lipschitz ODEs. In this sense, the different conditions of Assumptions \ref{Assumption2}-\ref{interconnection2} are natural extensions to study stability and/or recurrence properties in the stochastic set-valued hybrid case.
\section{Main Results}
\label{sec_results}
In this section, we present the main results of the paper. We will consider composite regular, locally Lipschitz certificate functions of the form
\begin{equation}\label{foster_function}
E_{\theta}(y):=(1-\theta) V(x)+\theta W(x,z), 
\end{equation}
where $\theta\in(0,1)$, and $V,W$ come from Section \ref{section3}. For convenience, we also introduce the set $L_{E_{\theta}}(c):=\{y\in\text{dom}~E_{\theta}:E_{\theta}(y)=c\}$, and the following quantities
\begin{equation}\label{epsiloncondition}
\varepsilon^*:= \frac{k_xk_z}{k_2k_z+k_1k_3},~~~\theta^*=\frac{k_3}{k_1+k_3},
\end{equation}
where the positive constants $(k_1,k_2,k_3,k_x,k_z)$ were introduced in Assumptions \ref{Assumption2}, \ref{Assumption3} and \ref{assumption6}. 
\subsection{Uniform Global Stability in Probability}
Given a compact set $\mathcal{A}\subset\mathbb{R}^{n_1}$, we are interested in establishing suitable stochastic stability properties for the SP-SHDS \eqref{SPSHDS1} with respect to the set
\begin{equation}\label{compactsetcomplete}
\tilde{\mathcal{A}}:=\left\{y\in\mathbb{R}^{n_1+n_2}:x\in\mathcal{A},~z\in \mathcal{M}(\tcb{x})\right\},
\end{equation}
which is compact due to Assumption \ref{manifold}. Among all the different probabilistic notions of stability, we will focus on the property of \emph{Uniform Global Stability in Probability}, introduced in \tcb{\cite[Sec. 2.3]{Teel_ANRC}}.

\vspace{0.1cm}
\begin{definition}\label{SHDS2def}
Consider the SP-SHDS \eqref{SPSHDS1} with state $y=(x,z)$, and a compact set $\mathcal{\tilde{A}}\subset\mathbb{R}^{n_1+n_2}$.
\begin{enumerate}[(C1)]
\item The set $\mathcal{\tilde{A}}$ is said to be \textsl{Uniformly Lyapunov stable in Probability}  if for each $\epsilon>0$ and $\rho>0$ there exists a $\delta>0$ such that for all initial conditions $\mathbf{y}_\omega(0,0) \in \mathcal{\tilde{A}}+\delta \mathbb{B}$, every maximal random solution ${\bf y}_\omega$ satisfies:
\begin{align}\label{stabilityprobability}
&\mathbb{P}\Big(\mathbf{y}_{\omega}(t,j)\in\mathcal{\tilde{A}}+\epsilon\mathbb{B}^\circ,~\forall~ (t,j)\in\text{dom}(\mathbf{y}_{\omega})\Big)\nonumber\\
&~~~~~~~~~~~~~~~~~~~~~~~~~~~~~~~~~~~~~~~~~~~~\geq 1-\rho.
\end{align}
\item The set $\mathcal{\tilde{A}}$ is said to be \textsl{Uniformly Lagrange stable in probability} if for each $\delta>0$ and $\rho >0$, there exists $\epsilon>0$ such that inequality \eqref{stabilityprobability} holds. 
\item The set $\mathcal{\tilde{A}}$ is said to be \textsl{Uniformly Globally Attractive in Probability} if for each $\epsilon>0, \rho>0$ and $R>0$, there exists $T\geq 0$ such that for all random solutions  ${\bf y}_\omega$ with $\mathbf{y}_\omega(0,0)\in \mathcal{\tilde{A}}+R \mathbb{B}$ the following holds: 
\begin{align*}
\mathbb{P}\Big(&\mathbf{y}_{\omega}(t,j)\in\mathcal{\tilde{A}}+\epsilon\mathbb{B}^\circ,\forall~t+j\geq T,(t,j)\in \text{dom}(\mathbf{y}_{\omega})\Big)\notag \\
&~~~~~~~~~~~~~~~~~~~~~~~~~~~~~~~~~~~~~~~~~~~~\geq 1-\rho.
\end{align*}
\end{enumerate}
If conditions (C1), (C2), and (C3) hold, system \eqref{SPSHDS1} is said to render the compact set $\mathcal{\tilde{A}}$ \textsl{Uniformly Globally Asymptotically Stable in Probability} (UGASp). \QEDB
\end{definition}

Our first result provides different sufficient conditions to guarantee UGASp of the set $\tilde{\mathcal{A}}$ for the SP-SHDS \eqref{SPSHDS1}. All the proofs are omitted due to space limitations.

\vspace{0.1cm}
\begin{theorem}\label{theorem1}
Let $\mathcal{A}\subset\mathbb{R}^{n_1}$ be compact, and $\tilde{\mathcal{A}}$ be given by \eqref{compactsetcomplete}. Suppose that Assumption \ref{manifold} holds, $\varepsilon\in(0,\varepsilon^*)$ and:
\begin{enumerate}[(a)]
\item Assumptions \ref{Assumption2}, \ref{Assumption3} and \ref{assumption6} hold with $\alpha_i\in\mathcal{K}_{\infty}$ for all $i\in\{1,2,3,4\}$, $\varphi_x\in\mathcal{P}\mathcal{D}(\mathcal{A})$, $\varphi_z\in\mathcal{P}\mathcal{D}$, and $\mu=0$.
\item There exists a function $\hat{\rho}\in\mathcal{P}\mathcal{D}(\tilde{\mathcal{A}})$ such that
\begin{equation}\label{jump_conditionproof1}
\int_{\mathbb{R}^m}\sup_{\substack{g\in G(y,v)}}~E_{\theta^*}(g)\mu(dv)\leq E_{\theta^*}(y)-\hat{\rho}(y),
\end{equation}
for all $y\in D$, where $E_{\theta^*}$ is given by \eqref{foster_function}.
\end{enumerate}
Then, system \eqref{SPSHDSreduced} renders the set $\tilde{\mathcal{A}}$ UGASp. 
\QEDB
\end{theorem}
The result of Theorem \ref{theorem1} relies on showing that $E_{\theta^*}$ is a
\emph{strong regular Lyapunov-Foster function} for the SP-SHDS \eqref{SPSHDS1} whenever $\varepsilon$ is sufficiently small. However, in some applications it might not be easy to find functions $V,W$ that satisfy all the conditions of Theorem \ref{theorem1} to construct a strong Lyapunov-Foster function. In that case, some of the assumptions can be relaxed if certain solutions of \eqref{SPSHDS1} can be ruled out. The next result addresses this case.

\vspace{0.1cm}
\begin{theorem}\label{theorem2}
Let $\mathcal{A}\subset\mathbb{R}^{n_1}$ be compact, and $\tilde{\mathcal{A}}$ be given by \eqref{compactsetcomplete}. Suppose that Assumption \ref{manifold} holds, $\varepsilon\in(0,\varepsilon^*)$, and 
\begin{enumerate}[(a)]
\item Assumptions \ref{Assumption2}, \ref{Assumption3}, and \ref{assumption6} hold with $\alpha_i\in\mathcal{K}_{\infty}$, for all $i\in\{1,2,3,4\}$, $\varphi_x\in\mathcal{P}_s\mathcal{D}(\mathcal{A})$, $\varphi_z\in\mathcal{P}_s\mathcal{D}$, and $\mu=0$.
\item At least one of the following conditions hold:
\begin{enumerate}[(1)]
\item Assumptions \ref{assumptionjumpsreduced} and \ref{jumpcoupledcondition} with $\mu=0$, $\rho_x=\rho_4$, and
\begin{equation}\label{relaxcondition1}
\frac{k_3k_4}{k_1}<c_x.
\end{equation}
\item Assumptions \ref{assumptionjumpsfast} and \ref{interconnection2} with $\mu=0$, $\rho_z=\rho_5$, and
\begin{equation}\label{relaxcondition2}
\frac{k_1k_5}{k_3}<c_z.
\end{equation}
\end{enumerate}
\item There does not exist an almost surely complete random solution $\bf y_\omega=(\bf x_\omega,\bf z_\omega)$ that remains in $L_{E_{\theta^*}}(c)$ for every $c>0$ for which $L_{E_{\theta^*}}(c)$ is non-empty.

\end{enumerate}
Then, system \eqref{SPSHDS1} renders the set $\tilde{\mathcal{A}}$ UGASp . \QEDB
\end{theorem}
\subsection{Uniform Global Recurrence}

\vspace{-0.2cm}
The second main property that we study in this paper is the property of uniform global recurrence, introduced in \tcb{\cite[Sec. 2.4]{Teel_ANRC}}. Compared to stability, this weaker property is commonly studied in stochastic systems for which a stable set might not exist. 
\vspace{0.1cm}
\begin{definition}\label{definition3}
An open, bounded set $\mathcal{O}\subset\mathbb{R}^{n_1+n_2}$ is \emph{Uniformly Globally Recurrent} (UGR) for the SP-SHDS \eqref{SPSHDS1} if there are no finite escape times for \eqref{flow_dynamics1} and for each $\rho>0$ and $R>0$ there exists $\tau\geq0$ such that every maximal random solution $\bf{y}_{\omega}=(\bf{x}_{\omega},\bf{z}_{\omega})$ starting in the set $R\mathbb{B}$ satisfies
\begin{align*}
&\mathbb{P}\Big(\big(\text{graph}({\bf y}_{\bf{\omega}})\subset (\Gamma_{<\tau}\times\mathbb{R}^n)\big)\lor \big(\text{graph}({\bf y}_{\bf\omega})\cap (\Gamma_{\leq\tau}\times\mathcal{O})\big)\Big)\\
&~~~~~~~~~~~~~~~~~~~~~~~~~~~~~~~~~~~~~~~~~~~~~~~~~~~~~~~~~~\geq 1-\rho.
\end{align*}
where $\text{graph}({\bf y_\omega}):=\big\{(t,j,s):s={\bf y}_\omega, (t,j)\in\text{dom}(y)\big\}$, and $\Gamma_{<\tau}:=\{(s,t)\in\mathbb{R}^2:s+t<\tau\}$.  \QEDB
\end{definition}
Loosely speaking, Definition \ref{definition3} says that from every initial condition, solutions to \eqref{SPSHDS1} either stop or hit the set $\mathcal{O}$, with a hitting time that is uniform over compact sets of initial conditions, and the solutions do not have finite escape times.

The property of UGR is studied mainly for open sets (relative to $C\cup D$) in order to guarantee suitable uniformity and robustness properties. However, since set-valued mappings satisfying the Basic Conditions (or even continuous functions) might not map open sets to open sets, we now consider SP-SHDS with quasi-steady state mappings characterized by continuous \emph{open functions}, which are continuous functions that map open subsets of their domain to open subsets of their codomain, see \tcb{\cite{OpenFunctionsCite}}. 

\begin{assumption}\label{open_function}
There exists a continuous open function $m:\mathbb{R}^{n_1}\to\mathbb{R}^{n_2}$ such that for all $x\in C_x$ we have that $h(x)\in C_z$. \QEDB
\end{assumption}
Under Assumption \eqref{open_function}, for each open and bounded set $\mathcal{O}\subset\mathbb{R}^{n_1}$, the set $$\tilde{\mathcal{O}}:=\{y\in\mathbb{R}^{n_1+n_2}:x\in\mathcal{O},z=h(x)\}$$ is also open and bounded. Using this fact, the next result parallels Theorem \ref{theorem1}, and establishes UGR of $\mathcal{O}$ for the SP-SHDS \eqref{SPSHDS1} when $\varepsilon$ is sufficiently small.

\vspace{0.1cm}
\begin{theorem}\label{theorem3}
Let $\mathcal{O}\subset\mathbb{R}^{n_1}$ be an open and bounded set. Suppose that Assumption \ref{open_function} holds, $\varepsilon\in(0,\varepsilon^*)$, and:
\begin{enumerate}[(a)]
\item Assumptions \ref{Assumption2}, \ref{Assumption3} and \ref{assumption6} hold with $\mathcal{M}=m$, $\mu_F>0$, and $\mathcal{A}=\overline{\mathcal{O}}$.
\item There exists a continuous function $\hat{\rho}:\mathbb{R}^{n_1+n_2}\to\mathbb{R}_{>0}$ and $\mu_J>0$ such that
\begin{equation*}
\int_{\mathbb{R}^m}\sup_{\substack{g\in G(y,v)}}~E_{\theta^*}(g)\mu(dv)\leq E_{\theta^*}(y)-\hat{\rho}(y)+\mu\mathbb{I}_{\tilde{\mathcal{O}}},
\end{equation*}
for all $y\in D$, where $E_{\theta^*}$ is given by \eqref{foster_function}.
\end{enumerate}

Then, system \eqref{SPSHDS1} renders the set $\tilde{\mathcal{O}}$ UGR. \QEDB
\end{theorem}
Similar to Theorem \ref{theorem2}, we can relax some of the strong decrease conditions on $E_{\theta^*}$ in order to establish UGR of $\tilde{\mathcal{O}}$, provided $\varepsilon$ is sufficiently small, and certain solutions can be ruled out from \eqref{SPSHDS1}. 

\vspace{0.1cm}
\begin{theorem}\label{theorem4}
Let $\mathcal{O}\subset\mathbb{R}^{n_1}$ be an open and bounded set. Suppose that Assumption \ref{open_function} holds, $\varepsilon\in(0,\varepsilon^*)$, and: 
\begin{enumerate}[(a)]
\item Assumptions \ref{Assumption2}, \ref{Assumption3}, and \ref{assumption6} hold with $\mathcal{M}=m$ and $\mathcal{A}=\overline{\mathcal{O}}$.
\item At least one of the following conditions is satisfied:
\begin{itemize}
\item Assumptions \ref{assumptionjumpsreduced} and \ref{jumpcoupledcondition} hold with $\rho_x=\rho_4$, and inequality \eqref{relaxcondition1} holds.
\item Assumptions \ref{assumptionjumpsfast} and \ref{interconnection2} hold with $\rho_z=\rho_5$, and inequality \eqref{relaxcondition2} holds.
\end{itemize}
\vspace{0.1cm}
\item There does not exist an almost surely complete random solution $\bf y_{\omega}=(\bf x_{\omega},\bf z_{\omega})$ that remains almost surely in the set $L_{E_{\theta^*}}(c)\cap (\mathbb{R}^{n_1+n_2}\backslash\tilde{\mathcal{O}})$ for every $c\geq0$ for which $L_{E_{\theta^*}}(c)\cap (\mathbb{R}^{n_1+n_2}\backslash\tilde{\mathcal{O}})$ is non-empty.
\end{enumerate}
Then, system \eqref{SPSHDS1} renders the set $\tilde{\mathcal{O}}$ UGR. \QEDB
\end{theorem}
%
\section{Conclusions}
\label{sec_conclusions}
We have introduced various Lyapunov-based conditions that can be used to certify the properties of uniform global asymptotic stability in probability and uniform global recurrence in a class of stochastic hybrid dynamical systems with multiple time scales in the flows. Our results rely on the construction of composite regular Lagrange-Foster and Lyapunov-Foster functions using simpler functions available for the reduced and boundary layer dynamics of the system. The tools that we have presented in this paper have the potential to enable the design and analysis of new types of multi-time scale feedback controllers and algorithms that can incorporate hybrid dynamics and stochasticity to overcome some of the limitations of traditional smooth approaches. Future work will focus on the synthesis of such algorithms.
\bibliography{dissertationA.bib}

\end{document}